\documentclass{amsart}

\usepackage{xcolor,amsmath,amsthm}
\usepackage[verbose]{hyperref}
\hypersetup{colorlinks=false,allbordercolors=blue,pdfborderstyle={/S/U/W 1}}

\theoremstyle{definition}

\newcommand{\LL}{\mathcal{L}}

\newcommand{\F}{\mathbf{F}}
\newcommand{\C}{\mathbf{C}}
\newcommand{\Q}{\mathbf{Q}}
\newcommand{\G}{\mathrm{G}}
\newcommand{\kG}{k[\mathrm{G}]}
\newcommand{\Z}{\mathbf{Z}}
\newcommand{\Ab}{\mathsf{Ab}}
\newcommand{\SL}{\operatorname{SL}}

\newcommand{\ddef}{\stackrel{{\rm \tiny{def}}}{=}}

 \input xy
 \xyoption{all}
 
\newtheorem{qu}{Question}
\newtheorem{exm}{Example}
\newtheorem{conj}{Conjecture}
\newtheorem{rmk}{Remark}

\begin{document}

\title{On the modular Plesken Lie algebra}

\author{John Cullinan}

\keywords{Lie algebra; modular representation; finite group.}

%{AMS Subject Classification: 20C15, 20C20}

\maketitle

\begin{abstract}
Let $\G$ be a finite group.  The Plesken Lie algebra $\LL[\G]$ is a subalgebra of the complex group algebra $\C[\G]$ and admits a direct-sum decomposition into simple Lie algebras based on the ordinary character theory of $\G$.  In this paper we review the known results on $\LL[G]$ and related Lie algebras, as well as introduce a conjecture on a characteristic $p$ analog $\LL_p[\G]$, with a focus on when $p$ divides the order of $\G$.
\end{abstract}

\section{Introduction}

The aim of this paper is to review known results on a certain Lie algebra attached to a finite group using ordinary representation theory, as well as to suggest new directions of research that incorporate modular representations.  Let us begin with some definitions.

Let $k$ be a field, $\G$ a finite group, and $k[\G]$ the group algebra of $\G$ over $k$.  The alternating bilinear map $[~,~]: k[\G] \times k[\G] \to k[\G]$ defined by 
\[
[x,y] = xy-yx
\]
endows $k[\G]$ with the structure of a Lie algebra.  For each $g \in \G$, define the element $\widehat{g} \ddef g-g^{-1} \in \kG$.  Following \cite{ct}, the \emph{Plesken Lie algebra} $\LL[\G]$ is the Lie subalgebra of $\kG$ linearly spanned over $k$ by the $\widehat{g}$.  

This raises a very natural question: \emph{what is the isomorphism type of $\LL[\G]$?} In particular, we can ask whether there is a natural relation between the Lie algebra structure of $\LL[\G]$ and the group structure of $\G$.   Since $\LL[\G]$ might not be a simple Lie algebra, we state our main question more formally as follows.

\begin{qu}
What are the isomorphism types (as Lie algebras over $k$) of the simple composition factors of $\LL[\G]$? 
\end{qu}

A natural starting point is when $k=\C$, since there is a well-understood classical theory of complex simple Lie algebras.  Indeed, in \cite{ct} the authors show that when $k = \C$, the Lie algebra $\LL[\G]$ is semisimple (owing to the fact that $\C[\G]$ is a direct sum of matrix algebras) and admits a direct sum decomposition into Lie subalgebras according to the degrees and Frobenius-Schur indicators of the complex irreducible representations of $\G$.  We review their main theorem and give an illustrative example in Section \ref{ordinary} below. 

Since the structure of $\LL[\G]$ is known when $k=\C$, we switch perspective and ask what happens for other fields $k$, in particular when the characteristic of $k$ divides the order of $\G$.  Here, there appear to be interesting connections with the modular representation theory of $\G$ which have yet to be studied.  After reviewing the current literature below, we propose new questions and conjectures for further research based on extensive calculation.

\section{History of the Problem} %Before turning to the Plesken Lie algebra in characteristic $p$, 

Before working with $\LL[\G]$ directly, let us first consider $k[\G]$ as a Lie algebra.  If $k = \C$, then it is well known that $\C[\G]$ decomposes as a direct sum of matrix algebras, which coincides with the decomposition into Lie algebras.  Suppose now that the characteristic of $k$ divides the order of $\G$ (so that we are in the \emph{modular} setting).  Then  Passi, Passman, and Sehgal determined in \cite{passman}  the exact conditions under which $k[G]$ is a nilpotent Lie algebra, and under which it is a solvable Lie algebra; they showed that the nilpotency/solvability can be determined purely from the group-theoretic properties of $\G$, foreshadowing the potential for deep connections between the group theory of $\G$ and the Lie algebra theory of $\LL[\G]$. 

Roughly a decade later, Smirnov and Zalesskii in \cite{zalesski} systematically studied the Lie-algebraic properties  of $\LL[\G]$; to the best of our knowledge this is the first work that concerns $\LL[\G]$ explicitly.  In that paper the authors show (among other things) that 
\begin{enumerate}
\item if $k$ has characteristic 0, then $\LL[\G]$ is solvable implies $\G$ is solvable, and
\item if $k$ has characteristic different than 2, then $\LL[\G]$ is nilpotent implies $\G^2$ (the subgroup of $\G$ generated by squares) is nilpotent.
\end{enumerate}
We remark that item (a) does not generalize to positive characteristic (as they point out in \cite[\S5.1]{zalesski}).  Additionally, item (b) appears to be the first structure theorem relating the Lie-algebraic properties of $\LL[\G]$ to the group theory of $\G$ in the modular setting (excluding characteristic 2).  

Returning to the case $k = \C$, in 2008 Cohen and Taylor determined the direct-sum decomposition of $\LL[\G]$ into complex simple Lie algebras of classical type \cite{ct} (we review the exact decomposition in Section \ref{ordinary}).  Subsequently, in \cite{marin}, Marin considered a more general setup, which is also briefly addressed in \cite{zalesski}, as follows.  Let $\alpha: \G \to \C^{\times}$ be a multiplicative character, and let $\LL_\alpha[\G]$ denote the span of the elements $g-\alpha(g)g^{-1}$.  Then $\LL_\alpha[\G]$ is also a Lie subalgebra of $k[\G]$ when endowed with the bracket operation.  When $\alpha = \mathbf{1}$, we have $\LL_{\mathbf{1}}[\G] = \LL[\G]$.  The main result of \cite{marin} is the determination of $\LL_\alpha[\G]$ as a direct sum of complex, simple Lie algebras of classical type and contextualizes the study of the Plesken algebra within the field of harmonic analysis.  More recently, Chaudhuri determined in 2020 in \cite{chaudhuri} the direct-sum decomposition of $\LL[\G]$ when $k$ is an algebraic extension of $\Q$.   And in 2023 Arjun and Romeo studied the Lie algebra  representation theory of $\LL[\G]$ over $\C$ \cite{aejm}.

It is evident that there is an established and significant body of work surrounding the Plesken Lie algebra in characteristic 0, relating the ordinary representation theory and group structure of $\G$ to the classical Lie algebra structure of $\LL[\G]$.  Before addressing the modular theory in Section \ref{modular}, we review the main theorem of \cite{ct} along with some new generalizations. For the remainder of the paper we use \textsf{Atlas} \cite{atlas} notation for finite simple groups.

\section{The Ordinary Plesken Lie Algebra} \label{ordinary}

Recall the classification of simple Lie algebras over $\C$ of classical type in Table \ref{classic_lie}. For small values of $n$, we have the following coincidences: 
$\mathfrak{B}_1 = \mathfrak{C}_1 =  \mathfrak{D}_1 = \mathfrak{A}_1$, $\mathfrak{B}_2 = \mathfrak{C}_2$, $\mathfrak{D}_2 = \mathfrak{A}_1 \times \mathfrak{A}_1$, and  $\mathfrak{D}_3 = \mathfrak{A}_3$.

\begin{table}[ht]
\begin{center}
\begin{tabular}{cccc}
Label & Isomorphsim Type & Dimension & $n$\\
\hline
$\mathfrak{A}_n$ & $\mathfrak{sl}_{n+1}$ & $(n+1)^2-1$ & $n \geq1$ \\
$\mathfrak{B}_n$ & $\mathfrak{o}_{2n+1}$ & $2n^2+n$ & $n \geq 2$ \\
$\mathfrak{C}_n$ & $\mathfrak{sp}_{2n}$ & $2n^2+n$ & $n \geq3$ \\
$\mathfrak{D}_n$ & $\mathfrak{o}_{2n}$ & $2n^2-n$ & $n \geq4$ 
\end{tabular}
\caption{Simple Lie Algebras of Classical Type}\label{classic_lie}
\end{center}
\end{table}

In \cite{ct} the authors show that if $t$ is the number of involutions of $\G$, then 
\begin{align} \label{dim}
\dim \LL[\G] = (|\G| - t - 1)/2.
\end{align}
The irreducible characters of $\G$ are partitioned according to their Frobenius-Schur indicator.  We denote these sets by $\mathcal{X}_1$, $\mathcal{X}_{-1}$, and $\mathcal{X}_{0}$ and  %, as in Section \ref{notation}.  
 %With all this notation set, we now 
recall the direct-sum decomposition proved in \cite[Thm.~5.1]{ct}:
\begin{align} \label{ctform}
\LL[\G] = \bigoplus_{\chi \in \mathcal{X}_1} \mathfrak{o}_{\chi(1)}  \bigoplus_{\chi \in \mathcal{X}_{-1}} \mathfrak{sp}_{\chi(1)}  \bigoplus_{\chi \in \mathcal{X}_0}~' \mathfrak{gl}_{\chi(1)},
\end{align}
where the prime signifies that there is just one summand $\mathfrak{gl}_{\chi(1)}$ for each pair $\lbrace \chi, \overline{\chi} \rbrace$ from $\mathcal{X}_0$.  The summands of $\mathfrak{gl}_n$ type are not simple, but contain $\mathfrak{sl}_n$ as a codimension-1 summand.  Let us now consider an illustrative example that shows how one can generalize (\ref{ctform}) to finite fields in ordinary characteristic.

\begin{exm} \label{main_ex}
Let $\G = {\rm L}_2(8)$, the simple group of order $504 = 2^3\cdot 3^2 \cdot 7$.  Using \textsf{Atlas} notation \cite[p.~6]{atlas}, the character table of $\G$ is given in Table \ref{L28}.

\begin{table}[ht]
\begin{center}
\begin{small}
\begin{tabular}{c|crrrrrrrrr}
&ind & 1A&2A&3A&7A&B*2&C*4&9A&B*2&C*4\\
\hline
$\chi_1$ &+&1&1&1&1&1&1&1&1&1\\
$\chi_2$ &+&7&$-1$&$-$2&0&0&0&1&1&1 \\
$\chi_3$  &+&7&$-$1&1&0&0&0&$-$\texttt{y9}&*2&*4 \\
$\chi_4$&+&7&$-$1&1&0&0&0&*4&$-$\texttt{y9}&*2 \\
$\chi_5$&+&7&$-$1&1&0&0&0&*2&*4&$-$\texttt{y9}\\
$\chi_6$&+&8&0&$-$1&1&1&1&$-$1&$-$1&$-$1\\
$\chi_7$&+&9&1&0&\texttt{y7}&*2&*4&0&0&0\\
$\chi_8$&+&9&1&0&*4&\texttt{y7}&*2&0&0&0\\
$\chi_9$&+&9&1&0&*2&*4&\texttt{y7}&0&0&0
\end{tabular}
\end{small}
\caption{Character Table of ${\rm L}_2(8)$}\label{L28}
\end{center}
\end{table}
\noindent Here, \texttt{y9} is a root of the polynomial $f_9(x) = x^3-3x-1$ and \texttt{y7} a root of $f_7(x) = x^3+x^2-2x-1$.  Both \texttt{y9} and \texttt{y7} are linear combinations of certain roots of unity:
\begin{align*}
\texttt{y9} &= -\zeta_9^5 - \zeta_9^4\\
\texttt{y7} &=-\zeta_7^4 + \zeta_7^3,
\end{align*}
and `*$m$' denotes replacing a root of unity $\zeta$ by $\zeta^m$. Let us now consider the Plesken Lie algebra of $\G$ from several perspectives.

\begin{enumerate}
\item[(1)] Applying (\ref{dim}) shows that $\dim \LL[\G]  = 220$.
\item[(2)] Applying (\ref{ctform}), we see that that $\LL[\G]$ admits the direct-sum decomposition
\[
\LL[\G]= \mathfrak{o}_7^4 \oplus \mathfrak{o}_8 \oplus \mathfrak{o}_9^3.
\]
\item[(3)] Let $p$ be a prime number coprime to $|G|$, so that $\F_p[\G]$ is semisimple (we call $p$ and \emph{ordinary} prime).  Define $\LL_p[\G]$ to be the Lie subalgebra of $\F_p[\G]$ spanned by the $\widehat{g}$.  Then the proof of \cite[Thm.~5.1]{ct} applies \emph{mutatis mutandis}, provided $\F_p$ contains enough roots of unity so that all irreducible representations of $\G$ over $\F_p$ are absolutely irreducible (we say that $\F_p$ is a \emph{splitting field} of $\G$).  If $\F_p$ is not a splitting field for $\G$, then the absolutely irreducible Galois-conjugate representations that are not realizable over $\F_p$ sum to an irreducible representation over $\F_p$.  Given this caveat, the proof of \cite[Thm.~5.1]{ct} carries through as well, though the dimensions of the simple factors are potentially different than over a splitting field.  We consider several special cases. \label{3}

\begin{enumerate}
\item \label{1st} Let $p = 71$.  Then $p$ is ordinary for $\G$ and both $f_9(x)$ and $f_7(x)$ split modulo $p$, and we have 
\[
\LL_p[\G]= \mathfrak{o}_7^4 \oplus \mathfrak{o}_8 \oplus \mathfrak{o}_9^3,
\]
where all Lie algebras in the sum are defined over the field $\F_{71}$.
\item \label{2nd} Let $p=17$.  Then $p$ is ordinary for $\G$, but only $f_9(x)$ splits modulo $p$, while $f_7(x)$ is irreducible.  The  three irreducible representations of degree 9 are not realizable over $\F_p$, but sum to an irreducible degree 27 representation.  Explicitly, this representation is given by the reduction modulo $17$ of the degree 27 representation over $\Z$ in the online \textsf{Atlas}:
\begin{center}
\href{https://brauer.maths.qmul.ac.uk/Atlas/v3/lin/L28/}{https://brauer.maths.qmul.ac.uk  /Atlas/v3/lin/L28/}
\end{center}
We then have
\[
\LL_p[\G]= \mathfrak{o}_7^4 \oplus \mathfrak{o}_8 \oplus \mathfrak{m},
\]
where $\mathfrak{m}$ is a 108-dimensional Lie algebra decomposing over $\F_{17}(\texttt{y7})$ into $\mathfrak{o}_9^{3}$.
\item \label{3rd} Let $p=5$.  Then $p$ is ordinary for $\G$, neither $f_9(x)$ nor $f_7(x)$ split modulo $p$, and we have 
\[
\LL_p[\G] = \mathfrak{o}_7  \oplus \mathfrak{o}_8 \oplus \mathfrak{l} \oplus \mathfrak{m},
\]
where $\mathfrak{l}$ is a 63-dimensional Lie algebra decomposing over $\F_5(\texttt{y9})$ into $\mathfrak{o}_7^{3}$ and $\mathfrak{m}$ is a 108-dimensional Lie algebra decomposing over $\F_5(\texttt{y7})$ into $\mathfrak{o}_9^{3}$.
\end{enumerate}
\end{enumerate}

We can generalize Items (\ref{1st}), (\ref{2nd}), and (\ref{3rd}) to number fields as follows.  Denote by  $K_7$ and $K_9$ the splitting fields over $\Q$ of $f_7(x)$ and $f_9(x)$, respectively.  Both $K_7$ and $K_9$ are cyclic extensions of $\Q$ of order 3 and the compositum $K_7K_9$ is Galois over $\Q$ with Galois group isomorphic to $C_3 \times C_3$.  

Let $\zeta_{63} \in \C$ be a primitive 63rd root of unity.  Then $\Q(\zeta_{63})$ is a splitting field of $\G$ over $\Q$, and contains the minimal splitting field $K_7K_9$.   Basic field theory shows that $[\Q(\zeta_{63}):\Q] = 36$ and a portion of the intermediate field diagram is as follows
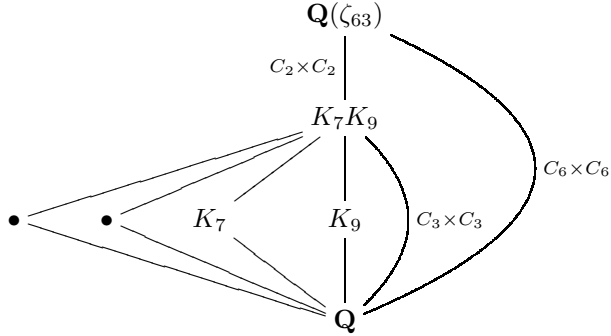
\[
\xymatrix{
&&&\Q(\zeta_{63}) \ar^{C_6 \times C_6}@/^6pc/@{-}[ddd] \ar_{C_2 \times C_2}@{-}[d]\\
&&&K_7K_9 \ar^{C_3 \times C_3}@/^2pc/@{-}[dd] \ar@{-}[dlll] \ar@{-}[dll] \ar@{-}[dl] \ar@{-}[d] \\
\bullet \ar@{-}[drrr]& \bullet \ar@{-}[drr] & K_7 \ar@{-}[dr] & K_9 \ar@{-}[d]\\
&&&\Q}
\]
where `$\bullet$' stands for a cubic extension of $\Q$ that does not concern this exposition (except for the calculations of the Chebatorev densities).

By the Chebatorev Density Theorem applied to the extension $K_7K_9/\Q$, it is a routine calculation to determine the proportion of primes for which the polynomials $f_7(x)$ and $f_9(x)$ split or are irreducible and, hence, the direct-sum decomposition of $\LL_p[\G]$. Specifically, we have

\begin{table}[ht]
\begin{center}
\begin{tabular}{cccc}
Proportion & Decomposition & $f_7(x)$ & $f_9(x)$ \\
of Primes & of $\LL_p[\G]$ && \\
\hline
1/9 & $\mathfrak{o}_7 \oplus \mathfrak{o}_8 \oplus \mathfrak{o}_7^3 \oplus \mathfrak{o}_9^3$ & split & split \\
4/9 & $ \mathfrak{o}_7  \oplus \mathfrak{o}_8 \oplus \mathfrak{l}_1\oplus \mathfrak{l}_2$ & irred. & irred. \\
2/9 & $\mathfrak{o}_7 \oplus \mathfrak{o}_8 \oplus \mathfrak{o}_7^3 \oplus \mathfrak{l}_2$ & split & irred. \\
2/9 & $\mathfrak{o}_7 \oplus \mathfrak{o}_8 \oplus \mathfrak{l}_1 \oplus \mathfrak{o}_9^3$ & irred & split
\end{tabular}
\caption{Splitting of $\LL_p[\G]$ for ordinary primes} \label{splitting_table}
\end{center}
\end{table}
\end{exm}

More generally, the minimal splitting field $K_\G$ of a finite group $\G$ is a subfield of the cyclotomic field $\Q(\zeta_m)$, where $m$ is the least common multiple of the orders of the elements of $\G$.  Thus, $K_\G$ is an abelian extension of $\Q$ and so class field theory would allow us to determine the densities of primes $p$ for which Galois-conjugate, absolutely irreducible representations of $\G$ are realizable modulo $p$, analogous to the case in Table \ref{splitting_table} above.  This is also a natural generalization of the results of \cite{chaudhuri} and \cite{ct} to characteristic $p$ fields when $p$ is ordinary.  

\section{The Modular Plesken Lie Algebra} \label{modular} Let $p$ be a prime dividing $|G|$ and let $k$ be an algebraically closed field of characteristic $p$.  Then $k[\G]$ is no longer semisimple, but decomposes as a direct sum of blocks.  The goal of this section is to introduce a conjecture on the structure of $\LL_p[\G]$ that takes into account the modular representation theory of $\G$ in a similar way to the ordinary case.  Observe that the dimension of $\LL_p[\G]$ is independent of whether $p$ is ordinary or modular.  

Let $\phi_1,\dots,\phi_s$ denote the Brauer characters of $\G$, where $s$ is the number of $p$-regular classes of $\G$.  As in the ordinary case, the set of Brauer characters is partitioned by Frobenius-Schur indicator into subsets $\mathcal{X}_1$, $\mathcal{X}_{-1}$, and $\mathcal{X}_0$.   The Lie algebra $\LL_p[\G]$ will generally not decompose into a direct sum of simple Lie algebras, but it does possess a composition series nonetheless.  Let us now set some notation in order to state a conjecture on the structure of $\LL_p[\G]$.  

Write $k[\G] = b_1 \oplus b_2 \oplus \cdots \oplus b_r$ as the block decomposition of $k[\G]$.  Then each $b_j$, being a two-sided ideal of $k[\G]$, assumes a Lie algebra structure over $k$ as well. Let us make some further remarks.

\begin{enumerate}
\item We write $\Ab(n)$ for an abelian Lie algebra of dimension $n$.  
\item If $p | (n+1)$, then the Lie algebra $\mathfrak{A}_n$ is not simple over $k$, since the scalar matrices have trace 0.  We write $\mathfrak{PA}_n$ for the codimension-1 simple factor of $\mathfrak{A}_n$ in this case.
\item In addition to the classical algebras, there exist exceptional simple Lie algebras in positive characteristic. When $p \geq 5$, the finite-dimensional simple Lie algebras over an algebraically closed field of characteristic $p$ have been fully classified.  Any such algebra is either of classical or Cartan type (Witt; Special; Hamiltonian; Contact) for $p \geq 7$ or, additionally, of Melikian type if $p=5$ \cite{ps}.  

Table \ref{eqtable} gives the dimensions of the exceptional  algebras and gives the conditions under which they are simple.With the exception of the Melikian algebras, these  Lie algebras are parameterized by $m,\underline{n}$ with $\underline{n} = [n_1,\dots,n_m] \in \Z_{>0}^m$.  Set $N = \sum n_i$. 

\begin{table}[ht]
\begin{center}
\begin{tabular}{|c|c|c|}
\hline
Lie algebra & Dimension & Simple \\
\hline
$W(m,\underline{n})$ &$p^N$& $p \ne 2$ and $m\ne 1$\\
\hline
$S(m,\underline{n})$  &$(m-1)p^N +1$& $m \geq 3$ \\
\hline
$H(m,\underline{n})$  &$p^N - 1$& $p >2$, $m \geq 2$\\
\hline
$K(m,\underline{n})$  &$p^N$ &$p>2$, $m \geq 3$\\
\hline
$M(n_1,n_2)$ &$5^{n_1 + n_2 + 1}$& $n_i >0$ \\
\hline
\end{tabular}
\caption{Exceptional Modular Lie Algebras}\label{eqtable}
\end{center}
\end{table}
In our Conjecture \ref{main_conj} below, the exceptional modular Lie algebras do not appear as simple factors. We conjecture that only the classical and abelian type Lie algebras appear as composition factors of $\LL[\G]$. 

\item We use lowercase fraktur letters for an arbitrary Lie algebra over $k$ and reserve $\mathfrak{A}_n$, $\mathfrak{B}_n$, $\mathfrak{C}_n$, $\mathfrak{D}_n$ for the classical Lie algebras of rank $n$.  For example, if $\LL_p[\G] = \mathfrak{a}_1 \oplus \mathfrak{a}_2 \oplus \mathfrak{A}_3$, then $\mathfrak{a}_1$ and $\mathfrak{a}_2$ are two unspecified Lie algebras over $k$, while $\mathfrak{A}_3$ refers to the Lie algebra of dimension 15 isomorphic to $\mathfrak{sl}_4$ over $k$.
\end{enumerate}

We are now ready to introduce our conjecture on the structure of $\LL_p[\G]$ in the modular setting.  The reader will notice that the degrees and indicators of the Brauer characters determine which classical Lie algebras appear in the decomposition, just as the the degrees and indicators of the complex characters did in the complex setting.  The main difference is that the degrees of the Brauer characters need not be the same as the complex characters and there are fewer Brauer characters than complex.  The remaining composition factors are all abelian, which is different from the complex setting.

\begin{conj} \label{main_conj}
Let $\G$ be a finite group and $p$ a modular prime for $\G$.  With all notation as above, we have 
\begin{enumerate}
\item $\LL_p[\G]$ admits a direct sum decomposition 
\[
\LL_p[\G] = \bigoplus_{j=1}^r \mathfrak{a}_j,
\]
where the $\mathfrak{a}_j$ denote the projection of $\LL_p[\G]$ onto the block $b_j$.  The Lie algebras $\mathfrak{a}_j$ are not necessarily simple.
\item The composition factors of $\LL_p[\G]$ are either abelian or of classical type.
\item \begin{enumerate}
\item for each $\phi \in \mathcal{X}_1$, there is a composition factor isomorphic to $\mathfrak{o}_{\phi(1)}$, and 
\item for each $\phi \in \mathcal{X}_{-1}$, there is a composition factor isomorphic to $\mathfrak{sp}_{\phi(1)}$, and 
\item for each pair $(\phi,\phi')$ of conjugate Brauer characters in $\mathcal{X}_0$, there is one composition factor isomorphic to $\mathfrak{gl}_{\phi(1)}$.  The algebra $\mathfrak{gl}_{\phi(1)}$ has a codimension-1 simple factor isomorphic to $\mathfrak{sl}_{\phi(1)}$, unless $p | \phi(1)$.  In that case, $\mathfrak{gl}_{\phi(1)}$ has a codimension-2 simple factor of type $\mathfrak{PA}_{\phi(1)-1}$. \label{brauer_decomp}
\end{enumerate}
\item A block of $k[\G]$ of defect 0 gives rise to direct summand of $\LL[\G]$ of classical type.
\end{enumerate}
\end{conj}

\begin{rmk}
The degrees of the Brauer characters do not necessarily agree with those of the ordinary characters, which can result in a very different overall Lie algebra structure of $\LL_p[\G]$ versus $\LL[\G]$, though the classical factors of $\LL[\G]$ as described in items (3a), (3b), (3c) of Conjecture \ref{main_conj} mimic the decomposition of (\ref{ctform}).
\end{rmk}

For the remainder of this section we give illustrative examples that are chosen to highlight different aspects of the conjecture.  In every case the reader will see that the simple composition factors can be predicted from the first two columns of the Brauer character table in the same manner as the ordinary $\LL[\G]$, as claimed in our conjecture. 

\begin{exm} \label{SL25ex}
Let $\G = \SL_2(5)$ and let $p=5$.  The indicator and degree columns of the Brauer and ordinary character tables are given in Table \ref{SL25table}.  
\begin{table}[ht]
\begin{center}
\begin{small}
\begin{tabular}{ccc||ccc}
Brauer&ind & degree&Ordinary & ind & degree \\
$\pmod 5$ && \\
\hline
$\phi_1$ & $+$ & 1 & $\chi_{1}$   &$+$&   1\\
$\phi_2$ & $-$ & 2 & $\chi_{2}$   &$-$&   2\\
$\phi_3$ & $+$ & 3 & $\chi_{3}$   &$-$&   2\\
$\phi_4$ & $-$ & 4 & $\chi_{4}$   &$+$&   3\\
$\phi_5$ & $+$ & 5 & $\chi_{5}$   &$+$&   3\\
&&&$\chi_{6}$   &$-$&   4\\
&&&$\chi_{7}$   &$+$&   4\\
&&&$\chi_{8}$   &$+$&   5\\
&&&$\chi_{9}$   &$-$&   6\\
\end{tabular}
\end{small}
\caption{Brauer -vs- Ordinary Characters of $\SL_2(5)$} \label{SL25table}
\end{center}
\end{table}
The group algebra $k[\G]$ decomposes into a direct sum of three blocks
\[
k[\G] = b_1 \oplus b_2 \oplus b_3,
\]
of defects 1, 1, and 0, respectively.  The ordinary characters $\chi_1$, $\chi_4$, $\chi_5$, and $\chi_7$ belong to $b_1$, while $\chi_2$, $\chi_3$, $\chi_6$, and $\chi_9$ belong to $b_2$.  The Steinberg character $\chi_8$ is the unique character belonging to $b_3$. 

In \textsf{Magma} we compute that $\LL_5[\G]$ admits the direct-sum decomposition
\[
\LL_5[\G] = \mathfrak{a}_1 \oplus \mathfrak{a_2} \oplus \mathfrak{B}_2,
\]
where $\mathfrak{a}_1$ and $\mathfrak{a}_2$ are Lie algebras with respective composition factors
\begin{align*}
\mathfrak{a}_1: &\ \Ab(3), \Ab(3), \Ab(3), \mathfrak{B}_1 \\
\mathfrak{a}_2: &\ \Ab(3), \Ab(3), \Ab(8), \Ab(10), \mathfrak{C}_1, \mathfrak{C}_2, 
\end{align*}
Observe that the simple composition factors $\mathfrak{B}_1$, $\mathfrak{C}_1$, $\mathfrak{B}_2$, $\mathfrak{C}_2$ that we computed in \textsf{Magma} agree part (3) of Conjecture \ref{main_conj}.  Observe also that the direct summand $\mathfrak{B}_2$ corresponds to the block $b_3$ of $k[\G]$ of defect 0.
\end{exm}

\begin{exm}
Let $\G = {\rm L}_2(8)$ as in Example \ref{main_ex}.  We select modular primes that are different from the defining characteristic (in contrast to Example \ref{SL25ex} above).  

\begin{table}[ht]
\begin{center}
\begin{small}
\begin{tabular}{ccc||ccc||ccc}
Brauer&ind & degree&Brauer & ind & degree & Ordinary & ind & degree \\
$\pmod 3$ &&& $\pmod{7}$&&&&& \\
\hline
$\phi_{1}$ & $+$ & 1 & $\phi_{1}$ & $+$ & 1 & $\chi_1$ & $+$ & 1 \\
$\phi_{2}$ & $+$ & 7 & $\phi_{2}$ & $+$ & 7 & $\chi_2$ & $+$ & 7 \\
$\phi_{3}$ & $+$ & 9 & $\phi_{3}$ & $+$ & 7 & $\chi_3$ & $+$ & 7 \\
$\phi_{4}$ & $+$ & 9 & $\phi_{4}$ & $+$ & 7 & $\chi_4$ & $+$ & 7 \\
$\phi_{5}$ & $+$ & 9 & $\phi_{5}$ & $+$ & 7 & $\chi_5$ & $+$ & 7 \\
&&& $\phi_{6}$ & $+$ & 8 & $\chi_6$ & $+$ & 8 \\
&&&&&& $\chi_7$ & $+$ & 9 \\
&&&&&&$\chi_8$ & $+$ & 9 \\
&&&&&& $\chi_9$ & $+$ & 9 
\end{tabular}
\end{small}
\caption{Brauer -vs- Ordinary Characters of ${\rm L}_2(8)$} \label{L28comp}
\end{center}
\end{table}

In characteristic 3, we have that $k[\G]$ decomposes into a sum of four blocks
\[
k[\G] = b_1 \oplus b_2 \oplus b_3 \oplus b_4,
\]
of defects 2, 0, 0, and 0, respectively.  The corresponding Lie algebra decomposition, according to \textsf{Magma}, is given by 
\[
\LL_3[\G] = \mathfrak{a}_1 \oplus \mathfrak{B}_4 \oplus \mathfrak{B}_4 \oplus \mathfrak{B}_4,
\]
where $\mathfrak{a}_1$ is a Lie algebra with composition factors
\[
\Ab(7), \Ab(21), \Ab(21), \Ab(21), \Ab(21), \mathfrak{B}_3.
\]
%Note that the Lie algebra $\mathfrak{B}_3$ is predicted from $\phi_2$ in the Brauer character table modulo 3. 

In characteristic 7, we have that $k[\G]$ decomposes into a sum of five blocks 
\[
k[\G] = b_1 \oplus b_2 \oplus b_3 \oplus b_4 \oplus b_5,
\]
of defects 1, 0, 0, 0, and 0, respectively.  The ordinary characters $\chi_1$, $\chi_6$, $\chi_7$, $\chi_8$, and $\chi_9$ belong to $b_1$, while $b_2$, $b_3$, $b_4$, $b_5$ contain only the characters $\chi_2$, $\chi_3$, $\chi_4$, $\chi_5$, respectively, since those characters remain irreducible modulo 7.  The corresponding Lie algebra decomposition, according to \textsf{Magma}, is given by
\[
\LL_7[\G] = \mathfrak{a}_1 \oplus \mathfrak{B}_3 \oplus \mathfrak{B}_3 \oplus \mathfrak{B}_3 \oplus \mathfrak{B}_3,
\]
where $\mathfrak{a}_1$ is a Lie algebra with composition factors 
\[
\Ab(8), \Ab(8), \Ab(8), \Ab(28), \Ab(28), \Ab(28), \mathfrak{D}_4.
\]
Again, observe that in both characteristic 3 and 7 the direct-sum decomposition and simple composition factors follow the predictions of parts (3) and (4) of Conjecture \ref{main_conj}.
\end{exm}

\begin{exm}
The purpose of this final example is to show that when there is a large power of $p$ dividing $\G$, the Brauer characters modulo $p$ may be quite different than the ordinary characters over $\C$.  %Nevertheless, the non-abelian composition factors as computed by \textsf{Magma} fit in with our conjectured structure theorem.
Let $\G = {\rm L}_2(25)$ and let $p=5$.  Then the dimensions and indicators of the Brauer and Ordinary characters are given in Table \ref{L225table}.  

\begin{table}[ht]
\begin{center}
\begin{small}
\begin{tabular}{ccc||ccc}
Brauer&ind & degree&Ordinary & ind & degree \\
$\pmod 5$ && \\
\hline
$\phi_{1}$& + &1 & $\chi_{1}$   &+&    1 \\
$\phi_{2}$& + &3 & $\chi_{2}$   &+&   13\\
$\phi_{3}$& + &3 & $\chi_{3}$   &+&   13\\
$\phi_{4}$& + &4 & $\chi_{4}$   &+&   24\\
$\phi_{5}$& + &5 & $\chi_{5}$   &+&   24\\
$\phi_{6}$& + &5 & $\chi_{6}$   &+&   24\\
$\phi_{7}$& + &8 & $\chi_{7}$   &+&   24\\
$\phi_{8}$& + &8 & $\chi_{8}$   &+&   24\\
$\phi_{9}$& + &9 & $\chi_{9}$   &+&   24\\
$\phi_{10}$& + &15 & $\chi_{10}$  &+&   25\\
$\phi_{11}$& + &15 & $\chi_{11}$  &+&   26\\
$\phi_{12}$& + &16 & $\chi_{12}$  &+&   26\\
$\phi_{13}$& + &25 & $\chi_{13}$  &+&   26\\
&&& $\chi_{14}$  &+&   26\\
&&& $\chi_{15}$  &+&   26
\end{tabular}
\end{small}
\caption{Brauer -vs- Ordinary Characters of ${\rm L}_2(25)$} \label{L225table}
\end{center}
\end{table}

The group algebra $k[\G]$ decomposes as a direct sum of two blocks $k[\G] = b_1 \oplus b_2$ of defect 2 and 0, respectively.  The Steinberg character $\chi_{10}$ belongs to $b_2$, while all others belong to $b_1$.  A calculation in \textsf{Magma} shows that the Lie algebra $\LL_5[\G]$ admits the decomposition 
\[
\LL_5[\G] = \mathfrak{a}_1 \oplus \mathfrak{B}_{12},
\]
where $\mathfrak{a}_1$ has composition factors as in Table \ref{SL225compfactors}.

\begin{table}[ht]
\begin{center}
\begin{small}
\begin{tabular}{|l|c|}
\hline
Composition Factor & Multiplicity \\
\hline
$\Ab(20), \Ab(25), \Ab(54), \Ab(90), \Ab(210)$  &1 \\
$\mathfrak{B}_4$, $\mathfrak{D}_8$ & \\
\hline 
$\Ab(16)$, $\Ab(24)$, $\Ab(30)$, $\Ab(48)$, $\Ab(64)$ & 2 \\
$\Ab(96)$, $\Ab(144)$, $\Ab(160)$, $\Ab(240)$ & \\
$\mathfrak{B}_2$, $\mathfrak{D}_4$, $\mathfrak{B}_7$ & \\
\hline
 $\Ab(56)$, $\Ab(120)$ & 3 \\
 \hline
 $\mathfrak{A}_1$ & 4 \\
\hline
$\Ab(6)$ & 7 \\
\hline
$\Ab(36)$ & 8 \\
\hline
\end{tabular}
\end{small}
\caption{Composition Factors of $\mathfrak{a}_1 \subseteq \LL_5[\G]$} \label{SL225compfactors}
\end{center}
\end{table}
\noindent Recall that $\mathfrak{B}_1 = \mathfrak{A}_1$ and $\mathfrak{D}_2 = \mathfrak{A}_1 \times \mathfrak{A}_1$, whence the four copies of $\mathfrak{A}_1$ in Table \ref{SL225compfactors}.
\end{exm}

We conclude the paper with further observations and questions for future researchers.

\begin{enumerate}
\item Observe that certain ordinary representations remain irreducible upon reduction modulo $p$.  These include the Steinberg representation, as well as any representation whose degree is a multiple of the largest power of $p$ dividing $|G|$. These representations occur in blocks of defect 0.  In the case of the Steinberg representation, the Frobenius-Schur indicator is $+1$, hence will give rise to a Lie algebra of type $\mathfrak{B}_n$ or $\mathfrak{D}_n$.

\item In addition to a proof of Conjecture \ref{main_conj}, it would be interesting to determine the precise number and dimensions of the abelian composition factors of $\LL_p[\G]$, as well as investigating any relation between the abelian factors and the defect groups of the blocks.
\end{enumerate}

\appendix

\section{Code and Data} All computations in this paper were performed in \textsf{Magma} \cite{magma}.  For the code below, ``\texttt{G}'' is a permutation group and ``\texttt{F}'' is a finite field.  This code will build the Plesken Lie Algebra ``\texttt{P}''. %, and from there the direct sum decomposition, composition series, etc.~can be found easily. \\

\noindent \texttt{
FG := GroupAlgebra(F, G);\\  						
L1,h:=Algebra(FG);  		\\					
FFG,f:=LieAlgebra(L1); 	\\							
L:=[]; 				\\						
for g in G do 			\\						
if Order(g) ne 2 and Order(g) ne 1 and Inverse(g) notin L then	\\
Append(\textasciitilde L,g); 									\\
end if; \\
end for;\\ 									
LL:=[]; 	\\								
for g in L do\\ 	
Append((\textasciitilde LL,FG!g-FG!Inverse(g) ));						\\	
end for; 									\\
LLL:=[]; 									\\
 for g in LL do 								\\	
Append((\textasciitilde LLL, (g@h)@f )); 					\\		
end for; 									\\
P:=sub<FFG|LLL>; 							
}


\begin{thebibliography}{99}

\bibitem{aejm} S.~Arjun, P.~Romeo. On representations of Plesken Lie algebras. Asian-Eur. J. Math. \textbf{16} (2023), no. 5, Paper No. 2350085.

\bibitem{chaudhuri} D.~Chaudhuri. Skew-symmetric elements of rational group algebras. Beitr. Algebra Geom. \textbf{61} (2020), no. 4, 719-729.

\bibitem{ct} A.~Cohen, D.E Taylor. On a Certain Lie Algebra Defined by a Finite Group. Amer. Math. Monthly \textbf{114} (2007), no. 7, 633-639.

\bibitem{atlas} J.~Conway \emph{et.~al.},~\textsl{ATLAS of Finite Groups}. Oxford University Press, Cambridge, 1985 

\bibitem{brauer_atlas} C.~Jansen \emph{et.~al.},~\textsl{An Atlas of Brauer Characters}, London Mathematical Society Monographs, Oxford University Press, New York, 1995

\bibitem{magma}
{W. Bosma}, {J. Cannon}, {C. Playoust}, \textit{The Magma algebra system.\ I.\ The user language}, J.\ Symbolic Comput.\ \textbf{24} (3--4), 1997, 235--265.

\bibitem{marin} I.~Marin. Group algebras of finite groups as Lie algebras. Comm. Algebra \textbf{38} (2010), no. 7, 2572-2584. 

\bibitem{passman} I.B.S.~Passi, D.~Passman, S.K.~Sehgal. Lie solvable group rings. Canad. J. Math. \textbf{25} (1973), 748-757. 

\bibitem{ps} A.~Premet, H.~Strade. Classification of finite dimensional simple Lie algebras in prime characteristics. \emph{Representations of algebraic groups, quantum groups, and Lie algebras}, 185-214, Contemp. Math., \textbf{413}, Amer. Math. Soc., Providence, RI, 2006.

\bibitem{zalesski}   M.B.~Smirnov, A.E.~Zalesskii. Lie algebra associated with linear group.
Comm. Algebra \textbf{9} (1981), no. 20, 2075-2100.

\end{thebibliography}
\end{document}